\documentclass[12pt]{amsart}
\usepackage{amssymb,amscd,amsthm,verbatim,nicefrac,amsmath,color,fancyhdr,mathrsfs,braket,graphicx, turnstile,hyperref,indentfirst,csquotes,url,enumitem,amsfonts, mathtools,lipsum}
\usepackage[letterpaper, left=2.5cm, right=2.5cm, top=2.5cm,
bottom=2.5cm,dvips]{geometry}
\hypersetup{
    colorlinks=true,
    linkcolor=blue,
    filecolor=magenta,      
    urlcolor=cyan,
    citecolor=blue,
    pdftitle={Minimising the Peak to $p$-average ratio in the vectorial case},
    pdfpagemode=FullScreen,
    }

\numberwithin{equation}{subsection}
\let\oldsection\section
\renewcommand{\section}{
  \renewcommand{\theequation}{\thesection.\arabic{equation}}
  \oldsection}
\let\oldsubsection\subsection
\renewcommand{\subsection}{
  \renewcommand{\theequation}{\thesubsection.\arabic{equation}}
  \oldsubsection}

\title[On the minimisation of the Peak--to--average ratio]{On the minimisation of the Peak--to--average ratio}

\author{Nikos Katzourakis}

\address{Department of Mathematics and Statistics, University of Reading, Whiteknights Campus, Pepper Lane, Reading RG6 6AX, UNITED KINGDOM}

\email{n.katzourakis@reading.ac.uk}

\makeatletter
\@namedef{subjclassname@2020}{\textup{2020} Mathematics Subject Classification}
\makeatother

\subjclass[2020]{Primary 35J47; 35J60; Secondary 35D30; 35A15}

\keywords{Crest Factor; Peak-to-average ratio; Calculus of variations in $\mathrm L^{\infty}$; higher order problems; generalised solutions; fully nonlinear equations; implicit PDEs; Baire Category method.}

\thanks{\!\!\!\!\!\!\!\texttt{The author has been partially financially supported through the EPSRC grant EP/X017109/1.}}

\def\R{\mathbb{R}}

\def\H{\mathrm{H}}
\def\E{\mathrm{E}}

\def\D{\mathrm{D}}
\def\Om{\Omega}

\def\X{\mathrm{X}}

\def\p{\mathrm{p}}

\def\d{\mathrm{d}}

\renewcommand\O{\mathcal{O}}

\newcommand{\av}{-\hspace{-13pt}\displaystyle\int}

\renewcommand{\d}{\mathrm{d}}

\newcommand{\mL}{\mathcal{L}}

\newcommand{\noi}{\noindent}
\newcommand{\ms}{\medskip}

\newcommand{\al}{\alpha}

\newcommand{\Ga}{\Gamma}

\newcommand{\la}{\lambda}
\newcommand{\La}{\Lambda}

\newcommand{\larrow}{\longrightarrow}
\newcommand{\ot}{\otimes}

\renewcommand{\L}{\mathrm L}

\renewcommand{\p}{\partial}

\newcommand{\sub}{\subseteq}

\newcommand{\by}{\times}

\newcommand{\ess}{\mathrm{ess}}

\newcommand{\beq}{\begin{equation}}

\newcommand{\eeq}{\end{equation}}

\newtheorem{theorem}{Theorem}

\newtheorem{proposition}[theorem]{Proposition}

\theoremstyle{definition}

\newtheorem{remark}[theorem]{Remark}

\begin{document}

\maketitle

\begin{abstract}
Let $\Omega \Subset \mathbb R^n$ and a continuous function $\H : \Omega\times \mathbb R^N \times\mathbb R^{N n} \times \cdots  \by \mathbb R^{N n^{k}}_s \longrightarrow \mathbb R$ be given, where $n,k,N \in \mathbb N$. For $p\in [1,\infty]$, we consider the
functional
\[
\ \ \ \ \ \mathrm E_p(u) := \big\| \mathrm H \big(\cdot,u,\mathrm D u, \ldots, \mathrm D^ku \big) \big\|_{\mathrm L^p(\Omega)},\ \ \ u\in \mathrm W^{k,p}(\Omega;\mathbb R^N).
\] 
In this note we are interested in the $\L^\infty$ variational problem
\[
\tag{$*$} \label{*}
\ \ \ \ \ \mathrm C_{\infty,p}(u_\infty)\, =\, \inf \Big\{\mathrm C_{\infty,p}(u) \ : \ u\in \mathrm W^{k,\infty}_\varphi(\Omega;\mathbb R^N), \ \mathrm E_1(u)\neq 0 \Big\},
\] 
where $\varphi\in \mathrm W^{k,\infty}(\Omega;\R^N)$ is given, $p$ is fixed, and
\[
\mathrm C_{\infty,p}(u)\, := \, \frac{\mathrm E_\infty(u)}{\mathrm E_p(u)} .
\] 
The variational problem \eqref{*} is ill-posed. $\mathrm C_{\infty,2}$ is known as the ``Crest factor" and arises as the ``peak--to--average ratio" problem in various applications, including eg.\ nuclear reactors and signal processing in sound engineering. We solve \eqref{*} by characterising the set of minimisers as the set of strong solutions to the eigenvalue Dirichlet problem for the fully nonlinear PDE
\[
\tag{$**$} \label{**}
\left\{
\ \
\begin{array}{ll}
\big| \mathrm H \big(\cdot,u,\mathrm D u, \ldots, \mathrm D^ku \big) \big|= \Lambda, & \text{ a.e.\ in }\Om,
\\
u = \varphi,  & \text{ on }\p\Omega,\\
\mathrm D u = \mathrm D \varphi, & \text{ on }\partial\Omega,\\
\ \ \  \ \vdots \ \ \ &\ \ \ \ \   \vdots \ \ \ \ \\
\mathrm D^{k-1}u = \mathrm D^{k-1}\varphi, & \text{ on }\partial\Omega.
\end{array}
\right.
\]
Under appropriate assumptions for $\H$, we show existence of infinitely-many solutions $(u,\Lambda) \in \mathrm W^{k,\infty}_\varphi(\Omega;\mathbb R^N) \times [\Lambda_*,\infty)$ to \eqref{**} for $\Lambda_*\geq0$, by utilising the Baire Category method for implicit PDEs. In the case of $k=1$ and $n=N$, these assumptions do not require quasiconvexity.
\end{abstract}


\smallskip

\section{Introduction}

Let $n,k,N \in \mathbb N$, $p\in [1,\infty]$, and let us also fix an open bounded domain $\Omega \Subset \mathbb R^n$. Given a mapping $u : \R^n \supseteq \Om \larrow \R^N$, the derivatives (of first, second, and $k$-th order) will be denoted by
\[
\left\{ \ \ \   \begin{split}
    \D u=\big(\D_i u_\al\big)_{i\in\{1,...,n\}}^{\al\in\{1,...,N\}}&\ :\ \ \Omega\larrow \R^{Nn},\\
    \D^2u=\left(\D^2_{ij}u_\al\right)_{i,j\in\{1,...,n\}}^{\al\in\{1,...,N\}}&\ :\ \ \Omega \larrow  \R^{Nn^2}_s,\\
  \vdots  \ \ \ \ \ \   \ \ \ \ \ \ & \ \ \ \ \ \ \  \ \ \ \vdots
     \\
    \D^k u=\left(\D^k_{i_1...i_k} u_\al\right)_{i_1,\ldots,i_k \in\{1,...,n\}}^{\al\in\{1,...,N\}}&\ : \ \ \Omega \larrow  \R^{Nn^k}_s,\ \ \ \ \ \ \ \ \ 
  \end{split}
  \right.
\]
and they are valued into their respective (symmetric) tensor spaces, defined as 
\[
  \R^{N n^k}_s:=\bigg\{{\bf X}\in \R^N \!\ot \underbrace{\R^n\otimes\cdots\otimes\R^n}_{k \text{ times}}\ :\ {\bf X}_{\al i_1\cdots i_k}= {\bf X}_{\al \sigma(i_1\cdots i_k)}, \sigma\text{ permutation on }(i_1,\ldots,i_k)\bigg\}.
\]
Finally, ${\bf X} \!:\! {\bf Y}$ denotes the Euclidean (Frobenius) inner product in $\R^{Nn^{k}}_s$. For brevity, we introduce the following compact notation for the ``$k$-th order jet" of $u$:
\[
\D^{[k]}u \,:=\, \Big(\cdot, u, \D u, \D^2u,...,\D^k u \Big) \ : \ \ \Om \,\larrow \, \Omega\times \mathbb R^N \times\mathbb R^{N n} \times \cdots  \by \mathbb R^{N n^{k}}_s. 
\]
Given a continuous function 
\[
\H \ :\ \ \Omega\times \Big(\mathbb R^N \times\mathbb R^{N n} \times \cdots  \by \mathbb R^{N (n-1)^{k}}_s \Big) \by \mathbb R^{N n^{k}}_s \longrightarrow \mathbb R, 
\]
we will symbolise the arguments of the supremand $\H$ throughout this note as $\H(x,\mathrm X, \textbf X)$. Let us consider the functional
\beq
\label{1.1}
\ \ \ \ \ \mathrm E_p(u) := \big\| \mathrm H (\D^{[k]}u ) \big\|_{\mathrm L^p(\Omega)},\ \ \ u\in \mathrm W^{k,p}(\Omega;\R^N).
\eeq
Here $\mathrm W^{k,p}(\Omega;\R^N)$ stands for the standard Sobolev space of $k$-times weakly differentiable maps $u : \R^n \supseteq \Om \larrow \R^N$, which lie, together with all their partial derivatives up to $k$-th order, in the space $\mathrm L^p(\Om;\R^N)$. For technical simplicity in the estimates, we will use the rescaled $\mathrm L^p$-norms, when $1 \leq p<\infty$:
\[
\| f \|_{\mathrm L^p(\Om)} := \left(\, \av_\Om |f|^p \, \mathrm d \mL^n \! \right)^{\!1/p}.
\]
Here $\mL^n$ denotes the Lebesgue measure on $\R^n$. Our general functional space, measure theory and PDE notation is generally standard, or otherwise self-explanatory, as e.g.\ in \cite{KV}. Let us define the functional
\beq
\label{1.2}
\ \  \ \ \ \ \ \mathrm C_{\infty,p}(u)\, := \, \frac{\mathrm E_\infty(u)}{\mathrm E_p(u)} ,\ \ \ \ u\in \mathrm W^{k,\infty}(\Omega;\R^N),\ \mathrm E_1(u)\neq 0.
\eeq
Since $\mathrm E_1(u)=0$ if and only if $\mathrm E_p(u)=0$ if and only if $\mathrm H (\D^{[k]}u )=0$ a.e.\ on $\Om$, \eqref{1.2} is well defined.

Herein we are interested in the $\L^\infty$ variational problem
\beq
\label{1.3}
\ \ \ \ \ \mathrm C_{\infty,p}(u_\infty)\, =\, \inf \Big\{\mathrm C_{\infty,p}(u) \ : \ u\in \mathrm W^{k,\infty}_\varphi(\Omega;\mathbb R^N), \ \mathrm E_1(u)\neq 0 \Big\},
\eeq
where $\varphi\in \mathrm W^{k,\infty}(\Omega;\R^N)$ is a given mapping, $p\in [1,\infty)$ is fixed, and $\mathrm W^{k,\infty}_\varphi(\Omega;\mathbb R^N)$ symbolises the standard affine Sobolev space of zero-trace mappings
\[
\mathrm W^{k,\infty}_\varphi(\Omega;\mathbb R^N) \,:=\,  \varphi+\mathrm W^{k,\infty}_0(\Omega;\mathbb R^N).
\]
In other words, we are interested in finding a minimiser $u_\infty\in \mathrm W^{k,\infty}_\varphi(\Omega;\mathbb R^N)$ which realises the infimum in \eqref{1.3}. The primary motivation to study the problem  \eqref{1.3} comes from applications. When $p\in\{1,2\}$ in  \eqref{1.2}, the functionals $\mathrm C_{\infty,1}$ and $\mathrm C_{\infty,2}$ are known as the ``Crest factor" and arises as the ``peak--to--average ratio" problem in numerous applications, including eg.\ the modelling of nuclear reactors, engineering, signal and sound processing, se e.g.\ \cite{DCP, FCCM, J, Le, P, PP}. Further, if $\H(x,\cdot,\cdot)$ is for all $x\in\Om$ positively homogeneous of some degree $K>0$, then \eqref{1.3} is equivalent  {to the} following constrained $\mathrm L^\infty$ variational problem (with integral constraint):
\beq
\label{1.4}
\ \ \ \ \ \mathrm E_{\infty}(u_\infty)\, =\, \inf \Big\{\mathrm E_{\infty}(u) \ : \ u\in \mathrm W^{k,\infty}_\varphi(\Omega;\mathbb R^N), \ \mathrm E_p(u)=1\Big\}.
\eeq

Interestingly, the variational problem \eqref{1.3} is ill-posed. On the one hand, the admissible class is not weakly* closed in $\mathrm W^{k,\infty}_\varphi(\Omega;\mathbb R^N)$ in any reasonable manner. On the other hand, the functional $C_{\infty,p}$ is not sequentially weakly* lower-semicontinuous. Therefore, a priori the infimum in \eqref{1.3} might not be attained by a mapping in the admissible class. This implies that standard variational approaches, including the Direct Method of the Calculus of Variations, cannot straightfowardly be applied to yield the existence of a minimiser. The problem \eqref{1.3} can be seen as a kind of generalised eigenvalue problem for supremal functionals in the context of the Calculus of Variations in $\mathrm L^\infty$, in the spirit of the works \cite{Clark-Katzourakis-2023, Clark-Katzourakis-2024, K}. Here the Crest factor \eqref{1.2} plays the role of a generalised Rayleigh quotient, but unlike the classical Rayleigh quotient, technical complications noted earlier arise because the denominator involves derivatives of the same order as the numerator.

In this note we demonstrate how we can easily resolve \eqref{1.3} without involving neither standard variational methods, nor any generalised methods for functionals which are not weakly lower semicontinuous, typically involving Young measures, relaxations and lower semi-continuous envelopes. Instead, our approach is inspired by some recent progress on higher order supremal variational problems made in \cite{Katzourakis-Moser-2019, Katzourakis-Moser-2023, Katzourakis-Moser-2024-LocalMinimisers}. More specifically, we prove that set of minimisers of \eqref{1.3} can be characterised as the set of strong a.e.\ solutions $(u,\Lambda) \in \mathrm W^{k,\infty}_\varphi(\Omega;\mathbb R^N) \times [\La_*,\infty)$ to the eigenvalue Dirichlet problem for the next fully nonlinear PDE of $k$-th order, above some ``critical" eigenvalue $\La_*\geq 0$:
\beq
\label{1.5}
\left\{ \ \ \
\begin{array}{ll}
\big| \H \big(\cdot,u,\mathrm D u, \ldots, \mathrm D^ku \big) \big|= \Lambda, & \text{ a.e.\ in }\Om,
\\
u = \varphi,  & \text{ on }\p\Omega,\\
\mathrm D u = \mathrm D \varphi, & \text{ on }\partial\Omega,\\
\ \ \ \ \vdots \ \ \ \ \ & \ \ \ \ \vdots \ \ \ \ \\
\mathrm D^{k-1}u = \mathrm D^{k-1}\varphi, & \text{ on }\partial\Omega.
\end{array}
\right.
\eeq
Given our notation for jets, we can rewrite the $(k-1)$-th order Dirichlet problem \eqref{1.5} in the following abbreviated form as:
\[
\left\{\ \ \ 
\begin{array}{ll}
 \big| \H \big(\mathrm D^{[k]}u \big)  \big| = \Lambda, & \text{ a.e.\ in }\Om, \ms
\\
\D^{[k-1]}u = \D^{[k-1]}\varphi ,  & \text{ on }\p\Omega.
\end{array}
\right.
\]
Our main result in this regard is Theorem \ref{theorem1}, which fully characterises the equivalence between these two classes of mappings. The idea of the proof is as follows: noting that H\"older's inequality implies $\mathrm C_{\infty,p}(u)\geq 1$, we show that the infimum in fact is equal to one if and only if $\big| \H \big(\mathrm D^{[k]}u \big)  \big|$ is essentially constant on $\Om$. Of course, the proof of equivalence still leaves open the question when \eqref{1.5} is actually solvable, and in fact solvability of this problem is a necessary condition in Theorem \ref{theorem1}. To this end, we present two different results in this regard. Both results rely essentially on the Dacorogna-Marcellini Baire Category method for implicit PDEs, for which the standard reference is \cite{D-M}. 

The first result regarding the existence of solutions for \eqref{1.5}, Proposition \ref{proposition1}, is essentially a direct adaptation of results from \cite{D-M} to our setting, providing sufficient conditions for the solvability of \eqref{1.5}, which unsurprisingly include quasiconvexity and coercivity for $\H$ with respect to the argument for the leading derivative term. The least possible eigenvalue $\La_*$ for solvability is determined by the value of the nonlinear operator applied to boundary condition $\varphi \in \mathrm W^{k,\infty}(\Omega;\mathbb R^N)$. However, for reasons relevant to approximation results involving quasiconvex functions, one needs additional convexity assumptions to have a completely general boundary condition as above. Otherwise, the problem is solvable for piecewise smooth boundary conditions (for details see \cite{D-M} and subsequently). 

Our second result regarding the existence of solutions for \eqref{1.5}, Proposition \ref{lem:existenceofsolutionstoH=C}, is a less straightforward application of the Baire Category method, and applies to the case of equal dimensions $n=N$ and for order $k=1$ in \eqref{1.5}. It relies on the solvability of the singular value problem, and interestingly does not require quasiconvexity of $\H$ with respect to the leading order term, only some other structural assumptions, it does however apply only to piecewise smooth boundary conditions $\varphi$. One could perhaps remove the assumption that $n=N$ by arguing as in the paper \cite{CKP}, but it is not clear whether this method can yield general boundary conditions in this case too as in \cite{CKP}. Anyhow, we refrain from delving into this task on this occasion, as the sufficient conditions provided by Proposition \ref{proposition1} are rather natural.  

We conclude this introduction by placing this problem into a wider context. Supremal variational problems were first considered in the 1960s by the pioneer of this field, Aronsson (see \cite{Aronsson-1965, Aronsson-1966, Aronsson-1967, Aronsson-1984}). For scalar first-order $\mathrm L^\infty$ variational problems (in which the supremand depends on the gradient of real-valued functions and perhaps lower-order terms), the theory is fairly well developed. Without any attempt to be exhaustive, we may refer to the following interesting relevant papers: \cite{AP, A-B, BJ, BDP, BK, CDP, KZ, MWZ, PWZ, PP, PZ, RZ, Ribeiro-Zappale-2024}. The vectorial first-order case is considerably more challenging, and largely still under development. Our work herein is inspired by recent progress in the vectorial and higher order case \cite{Clark-Katzourakis-2024, Katzourakis-Moser-2023, Katzourakis-Moser-2024-LocalMinimisers}, in which relevant phenomena were exhibited.

Except for the intrinsic theoretical interest, one of the advantages of minimising supremal energies is that they provide a uniformly small pointwise energy, whereas minimising an integral energy may allow for large spikes of the maximum pointwise energy, even though the area under the graph could still be small. This approach could provide better models in applications where this difference is relevant. However, the development of the theory is extremely challenging, as suprema do not enjoy the measure-theoretic properties of integrals, and practically no standard method of integral functionals applies to supremal functionals. In particular, they are not Gateaux differentiable regardless of the regularity of the supremand $\H$, and global minimisers may not minimise on subdomains with respect to their own boundary conditions.

\section{Main results and proofs}

We begin by establishing two existence results for strong solutions to the fully nonlinear eigenvalue problem \eqref{1.5}. Let us recall the following concepts taken from \cite{D-M}, for the convenience of the reader:

\begin{itemize}

\item Let $\Omega \Subset \mathbb R^n$ be a bounded open set. A map $u \in \mathrm W^{k,\infty}({\Om};\R^N)$ is called a piecewise-$\mathrm C^k$ map and we write 
\[
u \in \mathrm C^k_{{\rm pw}}(\overline{\Om};\R^N),
\] 
if there exists a union of countably many disjoint open set $(\Om_i)_{i=1}^\infty \sub \Om$ such that $\mL^n\big(\Om\setminus (\cup_{i=1}^\infty \Om_i))=0$, namely whose complement in $\Om$ is a Lebesgue nullset, and $u|_{\Om_i} \in \mathrm C^k(\overline{\Om}_i;\R^N)$, namely $u$ extends to a map which is $\mathrm C^k$ up to the boundary of $\Om_i$.

\ms

\item A continuous function $F : \mathbb R^{N n^{k}}_s  \larrow \R$ is called (Morrey) quasiconvex if
\[
f(\textbf X) \, \leq \, \av_{[0,1]^n}  f\big(\textbf X + \D^k \psi\big) \, \d \mL^n,
\] 
$\text{ for all } \psi \in  \mathrm W^{k,\infty}_0({\Om};\R^N) \text{ and all }\textbf X \in \R^{N n^{k}}_s.$
\end{itemize}

We first have the next general result regarding the existence of solutions to the fully nonlinear eigenvalue Dirichlet problem \eqref{1.5} for arbitrary order and dimensions. It is based on a straightforward application of results from \cite{D-M} on the Dacorogna-Marcellini Baire Category method for implicit PDE and differential inclusions.

\begin{proposition} \label{proposition1} Let $\Omega \Subset \mathbb R^n$ be a bounded open set, $n,k,N \in \mathbb N$ and 
\[
\H \ :\ \  \Omega\times \Big(\mathbb R^N \times\mathbb R^{N n} \times \cdots \by \mathbb R^{N n^{k-1}}_s\Big) \by \mathbb R^{N n^{k}}_s \longrightarrow \mathbb R
\]
a continuous function. Fix $\varphi \in  \mathrm W^{k,\infty}(\Omega;\R^N)$ and consider the eigenvalue Dirichlet problem \eqref{1.5}, namely we are looking for pairs $(u,\Lambda) \in \mathrm W^{k,\infty}(\Omega;\mathbb R^N) \times [0,\infty)$ such that
\beq
\label{2.1}
\left\{\ \ \ 
\begin{array}{rr}
 \big| \H \big(\mathrm D^{[k]}u \big)  \big| = \Lambda, & \text{a.e.\ in }\Om, \ms
\\
\D^{[k-1]}u = \D^{[k-1]}\varphi ,  & \text{ on }\p\Omega.
\end{array}
\right.
\eeq
Let us set
\beq \label{2.2}
 \La_* \, :=\, \underset{\Om}{\ess\,\sup}\,\big| \H \big(\mathrm D^{[k]}\varphi \big)  \big| ,
\eeq
that is, we take $\La_* = \E_\infty(\varphi)$.

Let us assume that: 
\beq
\label{2.3}
\left\{ \ \ \ 
\begin{split}
&\text{For any $(x,X) \in \Omega\times \mathbb R^N \times\mathbb R^{N n} \times \cdots  \by \mathbb R^{N n^{(k-1)}}_s$,}
\\
&\ \ |\H(x,X,\cdot)| \text{ is (Morrey) quasiconvex on }\R^{Nn^k}_s.
\end{split}
\right.
\eeq
Also,
\beq
\label{2.4}
\left\{ \ \ \ 
 {\begin{split}
& \ \text{There exists }\La\geq \La_* \text{ such that } \big\{ |\H(x,X,\cdot)| \leq \La\} \text{ is bounded in $\R^{Nn^k}_s$, for}
\\
&\text{\ any $x \in \Omega$ and any $X$ in a bounded subset of $\mathbb R^N \times\mathbb R^{N n} \times \cdots  \by \mathbb R^{N n^{(k-1)}}_s$.}
\end{split}}
\right.
\eeq
Then, the following are true:

\ms

\begin{enumerate}

\item[\emph{(1)}] If additionally $\varphi \in \mathrm C^k_{{\rm pw}}(\overline{\Om};\R^N)$ is piecewise-$\mathrm C^k$, then \eqref{2.1} has infinitely many solutions $u \in \mathrm W^{k,\infty}(\Omega;\mathbb R^N)$, where $\La$ as in \eqref{2.4}.

\ms

\item[\emph{(2)}]  If additionally for any $(x,X) \in  \Omega\times \big(\mathbb R^N \times\mathbb R^{N n} \times \cdots \by \mathbb R^{N n^{k-1}}_s\big)$, $|\H(x,X,\cdot)|$ is convex on $\R^{Nn^k}_s$, and $\La > \La_*$ in \eqref{2.4}, then \eqref{2.1} has infinitely many solutions $u \in \mathrm W^{k,\infty}(\Omega;\mathbb R^N)$, where $\La$ as in \eqref{2.4}.
\end{enumerate}

\ms

Suppose further
\beq
\label{2.5}
\left\{ \ \ \ 
 {\begin{split}
& \ \text{For all }\La\geq \La_*, \text{ the set } \big\{ |\H(x,X,\cdot)| \leq \La\} \text{ is bounded in $\R^{Nn^k}_s$, for any}
\\
&\text{ $x \in \Omega$ and any $X$ in a bounded subset of $\mathbb R^N \times\mathbb R^{N n} \times \cdots  \by \mathbb R^{N n^{(k-1)}}_s$.}
\end{split}}
\right.
\eeq

Then, the following are true:

\ms

\begin{enumerate}

\item[\emph{(3)}] If additionally $\varphi \in \mathrm C^k_{{\rm pw}}(\overline{\Om};\R^N)$ is piecewise-$\mathrm C^k$, then \eqref{2.1} has infinitely many solutions 
\[
\ \ \ \ (u,\La) \, \in \, \mathrm W_\varphi^{k,\infty}(\Omega;\mathbb R^N) \by [\La_*,\infty),
\]
and in fact, for any fixed $\La \geq \La_*$, it has infinitely many solutions $u \in \mathrm W_\varphi^{k,\infty}(\Omega;\mathbb R^N)$.

\item[\emph{(4)}] If additionally for any $(x,X) \in \Omega\times \mathbb R^N \times\mathbb R^{N n} \times \cdots  \by \mathbb R^{N n^{(k-1)}}_s$, $|\H(x,X,\cdot)|$ is convex on $\R^{Nn^k}_s$, then \eqref{2.1} has infinitely many solutions 
\[
\ \ \ \ (u,\La) \, \in \, \mathrm W_\varphi^{k,\infty}(\Omega;\mathbb R^N) \by (\La_*,\infty),
\]
and in fact, for any fixed $\La > \La_*$, it has infinitely many solutions $u \in \mathrm W_\varphi^{k,\infty}(\Omega;\mathbb R^N)$.
\end{enumerate}

\end{proposition}

\ms

\noi \text{\bf Proof of Proposition \ref{proposition1}.} Given $\La_*$ as in \eqref{2.2}, we fix $\La \geq \La_*$ and define
\[
F \ :\ \  \Omega\times \Big(\mathbb R^N \times\mathbb R^{N n} \times \cdots \by \mathbb R^{N n^{k-1}}_s\Big) \by \mathbb R^{N n^{k}}_s  \larrow \R
\]
by setting
\[
F_\La(x,X,\textbf X) \, :=\, |\H(x,X,\textbf X)| - \La.
\]
It follows that $F$ is a continuous function, and by assumption \eqref{2.3}, for any 
\[
(x,X) \in \Omega\times \mathbb R^N \times\mathbb R^{N n} \times \cdots  \by \mathbb R^{N n^{(k-1)}}_s, 
\]
the function $F_\La(x,X,\cdot)$ is (Morrey) quasiconvex on $\R^{N n^{k}}_s$. Further, by assumption \eqref{2.4}, it follows that $\{ F_\La(x,X,\cdot) \leq 0\}$ is a bounded set in $\R^{N n^{k}}_s$, for all $x\in\Om$ and for all $X$ in any bounded subset of $\mathbb R^N \times\mathbb R^{N n} \times \cdots  \by \mathbb R^{N n^{(k-1)}}_s$. Further, by the definition \eqref{2.2} and the expression defining $F_\La$, it follows that
\[
F_\La(\D^{[k]}\varphi) \leq 0,\ \text{ a.e.\ in }\Om,
\]
because $\La \geq \La_*$.

\ms

\noi (1) Since by assumption $\varphi \in \mathrm C^k_{{\rm pw}}(\overline{\Om};\R^N)$, the conclusion is a consequence of \cite[Theorem 6.23, p.\ 164]{D-M} applied to $F_\La$ and $\varphi$ above. Indeed, by the above result for implicit PDEs, it follows that there exist infinitely many strong solutions $u\in  \varphi + \mathrm W_0^{k,\infty}(\Omega;\mathbb R^N)$ to the PDE
\[
F_\La(\D^{[k]}u) = 0,\ \text{ a.e.\ in }\Om.
\]
Hence, the PDE of \eqref{2.1} is satisfied in the strong sense
\[
|\H(\D^{[k]}u)| = \La,\ \text{ a.e.\ in }\Om,
\]
whilst the Dirichlet boundary condition $\D^j u = \D^j \varphi$ of \eqref{2.1} is also satisfied pointwise on $\p\Om$, for all indices $j\in \{0,...,k-1\}$. 
\ms

\noi (2) Since by assumption $|\H(x,X,\cdot)|$ is convex on $\R^{N n^{k}}_s$, it follows that the same is true for $F_\La(x,X,\cdot)$. Further, since we also have assumed that $\La > \La_*$, the conclusion is a consequence of \cite[Theorem 6.23, p.\ 164]{D-M}, combined with the approximation result for convex functions \cite[Theorem 10.19, p.\ 245]{D-M} applied to $F_\La$ and $\varphi$ above. 

\ms

\noi (3) It is the same argument as in item (1), with the exception that, in view of the stronger assumption \eqref{2.5}, the conclusion is true for all $\La \geq \La_*$.

\ms

\noi (4) It is the same argument as in item (2), with the exception that, in view of the stronger assumption \eqref{2.5}, the conclusion is true for all $\La > \La_*$.      \qed
\ms


In the important special case that the order of differentiability $k$ equals one and $N=n$, we have the next result which interestingly \emph{does not require quasiconvexity as a key assumption}. It is based on the solvability of the Dirichlet singular value problem via the Dacorogna-Marcellini Baire Category method, as expounded in \cite{D-M}, which serves as an auxiliary problem.

\begin{proposition}\label{lem:existenceofsolutionstoH=C}
Fix $\Om \Subset \R^n$ and $\varphi \in \mathrm C^1_{{\rm pw}}(\overline{\Om};\R^n)$. Suppose that $\H:\Omega\times\R^n\times\R^{n^2} \larrow\R$ satisfies
$$
\H(x,\X, {\bf X})=h\big(x,\X, {\bf X}^\top{\bf X} \big),
$$ 
for some $h\in C\big(\Omega\times\R^n\times\R^{n^2}_s\big)$. Assume also that for any $ (x,\X) \in\Omega\times\R^n$, the function $h(x,X,\cdot)$ is strictly increasing along the direction of the identity, namely $t\mapsto h\big(x,X,t\hspace{1pt}\mathbb{I}_n\big)$ is strictly increasing on $\R$. Further, we suppose that there exists $\al_0>0$ such that 
\begin{equation}
\label{2.7}
\underset {\Omega \times \R^n}{\sup} h \big( \cdot,\cdot,\al_0 \mathbb{I}_n\big) \, <\, \infty,
\end{equation}
and also that the following uniform coercivity assumption is satisfied
\begin{equation}
\label{2.7A}
\lim_{t\to\infty}\Big(\, \underset {\Omega \times \R^n}{\inf} h \big( \cdot,\cdot, t \hspace{1pt}\mathbb{I}_n\big) \Big) \, =\, \infty.
\end{equation}
Then, the eigenvalue Dirichlet problem
\begin{equation}
\label{2.9}
 {  \begin{cases}
\ \  \H(\cdot, u,\D u)=\La, \ & \text{a.e.\ on }\Omega, \\
 \ \ u=\varphi,   \ & \text{on }\partial\Omega,
  \end{cases}}
\end{equation}
has (infinitely many) strong solutions $(u,\La)\in \mathrm W^{1,\infty}_\varphi(\Omega;\R^n) \by (\La_*,\infty)$, where
\begin{equation}
\label{2.8}
\La_* \, := \, \max\left\{ 
\underset {\Omega }{\mathrm{ess}\sup}\, h\Big(\cdot, \varphi,    \|\D \varphi\|^2_{\mathrm L^\infty(\Omega)} \mathbb{I}_n\Big) 
\, ,\, 
\underset {\Omega \times \R^n}{\sup} h \big( \cdot,\cdot,\al_0 \mathbb{I}_n\big)
 \right\}.
\end{equation}
In fact, for any fixed $\La > \La_*$, \eqref{2.9} has infinitely many solutions $u \in \mathrm W_\varphi^{k,\infty}(\Omega;\mathbb R^n)$. 

\smallskip

\noi Finally, the result above remains true if $h\geq0$, but we only require that $t\mapsto h\big(x,X,t\hspace{1pt}\mathbb{I}_n\big)$ be strictly increasing on $[0,\infty)$ (alongside \eqref{2.7} and \eqref{2.7A}), resulting in the exact same conclusions via the same proof, with the additional premise that $\La_*\geq0$. 
\end{proposition}

\noi {\bf Proof of Proposition }\ref{lem:existenceofsolutionstoH=C}. {\it \underline{Step 1}.} We first obtain some estimates directly from our assumptions. Let $\La_*$ be as in \eqref{2.8}, and fix $\La > \La_*$. It follows that
\[
\La \,>\, h\Big(x, \varphi(x),  \|\D \varphi\|^2_{\mathrm L^\infty(\Omega)} \mathbb{I}_n\Big),  \ \ \text{ a.e. }x\in \Omega.
\]
By the strict monotonicity assumption on $h$, we have
\[
\begin{split}
h\big(x, \varphi(x),(\cdot)\mathbb{I}_n\big)^{-1} (\La) &>\, \lVert\D \varphi \rVert_{\mathrm L^\infty(\Omega)}^2 \\
  &=\, \underset{\Om}{\ess\sup} \,\big\{ (\D \varphi^\top \D \varphi ) : \mathbb I_n\big\}
  \\
  & \geq \, \underset{\Om}{\ess\sup} \,\Big\{  \max_{e\in \R^n, |e|=1} (\D \varphi^\top \D \varphi ) : e\ot e \Big\} 
  \\ 
  &=\,  \underset{\Om}{\ess\sup} \, \la_n(\D \varphi^\top \D \varphi ),
  \end{split}
\]
pointwise on $\Om$, which implies
\beq
\label{2.11A}
h\big(x, \varphi(x),(\cdot)\mathbb{I}_n\big)^{-1} (\La)\, >\, \la_n(\D \varphi^\top \D \varphi ),\ \ \text{ a.e.\ on $\Omega$.}
\eeq
Further, again by \eqref{2.8}, we have
\[
\La \, > \, h \big( x,X, \al_0 \mathbb{I}_n\big),
\]
for some $\al_0>0$ and all $(x,X) \in \Omega \times \R^n$. By the monotonicity assumption on $h$ and the above inequality, we infer that
\[
\begin{split}
 h\big(x,X,(\cdot)\mathbb{I}_n\big)^{-1}(\La) \, >\, h\big(x,X,(\cdot)\mathbb{I}_n\big)^{-1}\Big(  h \big( x,X, \al_0 \mathbb{I}_n\big)\Big) \, =\, \al_0,
\end{split}
\]
for all $(x,X)\in\Om\by \R^n$, which yields
\beq
\label{2.12A}
\begin{split}
\underset{(x,X)\in\Om\by \R^n }{\inf} h\big(x,X,(\cdot)\mathbb{I}_n\big)^{-1}(\La) \, >\, \al_0.
\end{split}
\eeq 
Further, by the uniform coercivity assumption \eqref{2.7A} on $h$, there exists an $A\in\R$ such that
\[
\La \, \leq\, \underset{(x,X)\in\Om\by \R^n }{\inf}h\big(x,X,A\hspace{1pt}\mathbb{I}_n\big).
\]
By the monotonicity of $h$, this estimate is equivalent to
\beq
\label{2.13A}
\underset{(x,X)\in\Om\by \R^n }{\sup}h\big(x,X,(\cdot)\mathbb{I}_n\big)^{-1}(\La) \, \leq \, A, \ \ \ A\in\R.
\eeq
\noi {\it \underline{Step 2}.} We now demonstrate the solvability of \eqref{2.9} by constructing an auxiliary problem. Then, the estimates of {\it Step 1} will be the sufficient conditions for existence of solution to the auxiliary problem. To this aim, let us define the function $\al : \Omega \times \R^n \larrow \R$ by setting
\beq
\label{Al}
  \al(x,X) \,:= \, h\big(x,X,(\cdot)\mathbb{I}_n\big)^{-1}(\La).
\eeq
For any matrix $\X\in \R^{n^{2}}_s$, let $\{\lambda_1(\X)$, ...\,,$\lambda_n(\X)\}$ denote its eigenvalues in increasing order. Consider the following Dirichlet problem involving singular values of the gradient matrix
\begin{equation}
\label{5.5}
 {  \ \ \ \begin{cases}
\ \ \lambda_i(\D u^\top \D u)\,=\,\al(\cdot, u), & \text{a.e. on }\Omega, \quad i=1,\ldots,n, \\
 \ \   u= \varphi, \,  & \text{on }\partial\Omega.
  \end{cases} }
\end{equation}
By \cite[Theorem 7.28, p.\ 199]{D-M}, the problem \eqref{5.5} has (infinitely many) strong solutions in $\mathrm W^{1,\infty}_\varphi(\Omega;\R^n)$, as long as $\al$ is continuous, bounded above, bounded below by a positive constant away from zero, and the boundary condition $\varphi \in \mathrm C^1_{{\rm pw}}(\overline{\Om};\R^n)$ is a ``strict subsolution", in the sense that
\[
 \lambda_n(\D \varphi^\top \D \varphi) \,<\, \al(\cdot, \varphi ),  \ \ \text{ a.e.\ on }\Omega.
\]
By \eqref{2.11A}, \eqref{2.12A}, \eqref{2.13A} and \eqref{Al}, it follows that all the necessary conditions for the solvability of the problem \eqref{5.5} are satisfied.

\ms

\noi {\it \underline{Step 3}.} We conclude by showing how the solvability of the auxiliary problem \eqref{5.5} implies the solvability of \eqref{2.9}. For any fixed $u \in \mathrm W^{1,\infty}_\varphi(\Omega;\R^n)$ solving \eqref{5.5} for a given $\La>\La_*$, by the spectral theorem we have that there exists a measurable map of orthogonal matrices $\O \in L^\infty\big(\Omega;\mathrm O(n,\R)\big)$ such that
\begin{equation}
\label{5.8}
  \begin{aligned}
  \D u^\top \D u\, &=\, \O\begin{bmatrix}\lambda_1(\D u^\top \D u)& &\mathbb{O}\\ & \ddots & \\\mathbb{O}& &\lambda_n(\D u^\top \D u)\end{bmatrix}\O^\top  
  \\ 
  &=\, \al(\cdot, u)\hspace{1pt}\O\hspace{1pt}\mathbb{I}_n\hspace{1pt}\O^\top  \\ 
  &= \, \al(\cdot, u)\hspace{1pt}\mathbb{I}_n,
  \end{aligned}
\end{equation}
a.e.\ on $\Omega$. We may then compute
\[
\begin{split}
\H\big( x,\, u(x), \, \D u(x) \big) & = \, h\Big( x,\, u(x), \, \D u^\top \D u(x) \Big)
\\
& = \, h\Big(  x, \, u(x), \,   \al\big(x, u(x)\big) \mathbb{I}_n \Big)
\\
& =  \,h\Big(  x, \, u(x),  \,  (\cdot)\mathbb{I}_n \Big) \big(\al\big(x, u(x)\big)\big)
\\
& = \,\Big[ h\Big(  x, \, u(x),  \,  (\cdot)\mathbb{I}_n \Big) \circ h\Big(  x, \, u(x),  \,  (\cdot)\mathbb{I}_n \Big)^{-1}\Big](\La)
\\
& =  \,\La,
\end{split}
\]
a.e.\ on $\Omega$. The conclusion ensues. \qed
\ms


Now we may state and prove the main result in this note, regarding the equivalence between minimisers of the Crest factor functional in \eqref{1.3}, and the class of solutions to the eigenvalue Dirichlet problem  \eqref{1.5}.

\begin{theorem} \label{theorem1} Let $\Omega \Subset \mathbb R^n$ be a bounded open set, $n,k,N \in \mathbb N$ and 
\[
\H \ :\ \  \Omega\times \mathbb R^N \times\mathbb R^{N n} \times \cdots  \by \mathbb R^{N n^{k}}_s \longrightarrow \mathbb R
\]
a given continuous function. Let $\mathrm C_{\infty,p}$ be the Crest factor functional given by \eqref{1.1}-\eqref{1.2}, and fix $\varphi \in  \mathrm W^{k,\infty}(\Omega;\R^N)$. Suppose further that the problem \eqref{1.5} is solvable\footnote{Propositions \ref{proposition1} and \ref{lem:existenceofsolutionstoH=C} provide sufficient conditions on $\H,\varphi,\La_*$ for this to be true.} for all $\La\geq \La_*$, for some  {$\La_* >0$}.

\noi Then, for any $ u_\infty \in  \mathrm W^{k,\infty}_\varphi(\Omega;\R^N)$, the following statements are equivalent:

\begin{enumerate}

\item[\emph{(1)}] $u_\infty$ is minimiser of the variational problem \eqref{1.3}, namely
\[
\mathrm C_{\infty,p}(u_\infty)\, =\, \inf \Big\{\mathrm C_{\infty,p}(u) \ : \ u\in \mathrm W^{k,\infty}_\varphi(\Omega;\mathbb R^N), \ \mathrm E_1(u)\neq 0 \Big\}.
\]

\item[\emph{(2)}]  There exists $\La >0 $, such that $u_\infty$ is a strong a.e.\ solution to the Dirichlet problem \eqref{1.5}, namely
\[
\left\{\ \ \ 
\begin{array}{ll}
 \big| \H \big(\mathrm D^{[k]}u \big)  \big| = \Lambda, & \text{ a.e.\ in }\Om, \ms
\\
\D^{[k-1]}u = \D^{[k-1]}\varphi ,  & \text{ on }\p\Omega.
\end{array}
\right.
\]
\end{enumerate}
\end{theorem}

\noi {\bf Proof of Theorem} \ref{theorem1}. (2) $\Rightarrow$ (1):  {By assumption, we have that $(u_\infty,\La) \in \mathrm W^{k,\infty}_\varphi(\Om;\R^N)\by (0,\infty)$ solves the problem \eqref{1.5}}. This implies that $\big| \H \big(\mathrm D^{[k]}u_\infty \big) \big| = \Lambda$ a.e.\ on $\Om$, and also that $u_\infty - \varphi \in \mathrm W^{k,\infty}_0(\Omega;\mathbb R^N)$, or equivalently 
\[
u_\infty \in \mathrm W^{k,\infty}_\varphi(\Omega;\mathbb R^N). 
\]
 {In particular, the satisfaction of the PDE implies
\[
\La = \underset{\Om}{\ess \sup}\, \big| \H \big(\mathrm D^{[k]}u_\infty \big) \big| =\mathrm E_{\infty}(u_\infty),
\]
and also
\[
\La =  \left( \,\av_\Om \big| \H \big(\mathrm D^{[k]}u_\infty \big) \big|^p \, \d \mL^n \!\right)^{\!1/p}=\mathrm E_p(u_\infty),
\]
for any $p\geq 1$.} Note that, since we are using the rescaled $L^p$ norms on $\Om$, H\"older's inequality implies that
\[
\E_p (u) \leq \E_\infty (u), \ \text{ for all }u \in \mathrm W^{k,\infty}_\varphi(\Omega;\mathbb R^N).
\]
Therefore, by \eqref{1.2},
\beq
\mathrm C_{\infty,p}(u) \, \geq \, 1, \ \text{ for all }u \in \mathrm W^{k,\infty}_\varphi(\Omega;\mathbb R^N), \text{ when } \E_1 (u) \neq 0.
\eeq
Additionally, it follows that
\[
\mathrm C_{\infty,p}(u_\infty) = \frac{\mathrm E_{\infty}(u_\infty)}{\mathrm E_{p}(u_\infty)} = \frac{\big\| \H \big(\mathrm D^{[k]}u_\infty \big) \big\|_{\mathrm L^\infty(\Om)} }{ \big\| \mathrm H \big(\mathrm D^{[k]}u_\infty \big) \big\|_{\mathrm L^p(\Om)} } = \frac{\La}{\La}=1.
\]
Hence, $u_\infty$ realises the infimum in \eqref{1.3}, as claimed.

\ms

\noi (1) $\Rightarrow$ (2): Suppose $u_\infty$ realises the infimum in \eqref{1.3}. This implies that $u_\infty$ is in the admissible class, therefore $u_\infty \in \mathrm W^{k,\infty}_\varphi(\Omega;\mathbb R^N)$, or equivalently 
\[
u_\infty -\varphi \in \mathrm W^{k,\infty}_0(\Omega;\mathbb R^N),
\]
which implies that the $(k-1)$-th order Dirichlet boundary condition of \eqref{1.5} is satisfied. Further, $\mathrm E_1(u_\infty)\neq 0$. Equivalently, $\H(\D^{[k]}u_\infty)\neq 0$ on a Borel subset of positive Lebesgue measure on $\Om$,  {which in turn gives that $\mathrm E_\infty(u_\infty)>0$}. By H\"older's inequality, we again have
\[
\E_p (u) \leq \E_\infty (u), \ \text{ for all }u \in \mathrm W^{k,\infty}_\varphi(\Omega;\mathbb R^N),
\]
which implies that
\beq
\mathrm C_{\infty,p}(u) \, \geq \, 1, \ \text{ for all }u \in \mathrm W^{k,\infty}_\varphi(\Omega;\mathbb R^N),\ \text{when }  {\mathrm E_1(u)\neq 0}.
\eeq
We will now demonstrate that the infimum (realised by the energy value $u_\infty$) is precisely equal to $1$, namely,
\[
\mathrm C_{\infty,p}(u_\infty) = 1.
\]
To this aim, fix any solution $(u_0,\La_0) \in \mathrm W^{k,\infty}_\varphi(\Om;\R^N)\by [\La_*,\infty)$ to \eqref{1.5}, guaranteed to exist for some  {$\La_* >0$} by our assumptions. The arguments used above in the implication ``(2) $\Rightarrow$ (1)" (with $u_0$ in the place of $u_\infty$) prove that $u_0$ is a global minimiser with energy value equal to $1$, namely 
\[
\mathrm C_{\infty,p}(u_0)=1. 
\]
Since $u_\infty$ is by assumption a global minimiser as well, it follows that it must be at the same energy level, since
\[
\begin{split}
\mathrm C_{\infty,p}(u_\infty)\, &=\, \inf \Big\{\mathrm C_{\infty,p}(u) \ : \ u\in \mathrm W^{k,\infty}_\varphi(\Omega;\mathbb R^N), \ \mathrm E_1(u)\neq 0 \Big\} 
\\
& \leq\, \mathrm C_{\infty,p}(u_0)
\\
&=1.
\end{split}
\]
We will now conclude by showing that $u_\infty$ as well must be a strong solution to the PDE system in \eqref{1.5}. Since $\mathrm C_{\infty,p}(u_\infty)=1$, by \eqref{1.2} this implies
\[
\E_{\infty}(u_\infty) \,=\, \E_{p}(u_\infty).
\]
Therefore, by \eqref{1.1}
\[
\E_{\infty}(u_\infty)^p \,=\, \big\| \H \big(\mathrm D^{[k]}u_\infty \big) \big\|_{\mathrm L^p(\Om)}^p, 
\]
which yields
\beq
\label{1.8}
\E_{\infty}(u_\infty)^p \,=\, \av_\Om \big| \H \big(\mathrm D^{[k]}u_\infty \big) \big|^p \, \d \mL^n.
\eeq
Since
\[
\E_{\infty}(u_\infty) \,=\, \underset{\Om}{\ess \, \sup}\, \big| \H \big(\mathrm D^{[k]}u_\infty \big) \big| ,
\]
it follows that
\[
\big|\H \big(\mathrm D^{[k]}u_\infty \big) \big|^p \,\leq \, \E_{\infty}(u_\infty)^p \ \text{ a.e.\ on }\Om,
\]
and in particular we have that
\beq
\label{1.9}
 {\mL^n \bigg( \Om \setminus  \Big[ \Big\{ \big|\H \big(\mathrm D^{[k]}u_\infty \big) \big| \,= \, \E_{\infty}(u_\infty)\Big\} \bigcup \Big\{ \big|\H \big(\mathrm D^{[k]}u_\infty \big) \big| \,<\,  \E_{\infty}(u_\infty)\Big\} \Big] \bigg)=0.}
\eeq
By \eqref{1.8}-\eqref{1.9}, we have
\[
\begin{split}
\E_{\infty}(u_\infty)^p \, =& \ \av_\Om \big| \H \big(\mathrm D^{[k]}u_\infty \big) \big|^p \, \d \mL^n 
\\
=& \
 \frac{1}{\mL^n(\Om)}\int_{\big\{|\H (\mathrm D^{[k]}u_\infty)| =  \E_{\infty}(u_\infty)\big\}} \big| \H \big(\mathrm D^{[k]}u_\infty \big) \big|^p \, \d \mL^n 
 \\
 & \ +
  \frac{1}{\mL^n(\Om)}\int_{\big\{|\H (\mathrm D^{[k]}u_\infty)| <  \E_{\infty}(u_\infty)\big\}} \big| \H \big(\mathrm D^{[k]}u_\infty \big) \big|^p \, \d \mL^n
  \\
  =& \
 \frac{\E_{\infty}(u_\infty)^p}{\mL^n(\Om)} \mL^n\Big(\Big\{|\H (\mathrm D^{[k]}u_\infty)| =  \E_{\infty}(u_\infty)\Big\}\Big)  
 \\
 & \ +
  \frac{1}{\mL^n(\Om)}\int_{\big\{|\H (\mathrm D^{[k]}u_\infty)| <  \E_{\infty}(u_\infty)\big\}} \big| \H \big(\mathrm D^{[k]}u_\infty \big) \big|^p \, \d \mL^n. 
  \end{split}
\]
Therefore,
\[
\begin{split}
\int_{\big\{|\H (\mathrm D^{[k]}u_\infty)| <  \E_{\infty}(u_\infty)\big\}} \big| \H \big(\mathrm D^{[k]}u_\infty \big) \big|^p \, \d \mL^n   \ =& \
\mL^n(\Om) \E_{\infty}(u_\infty)^p 
\\
& \ - \E_{\infty}(u_\infty)^p \mL^n\Big(\Big\{|\H (\mathrm D^{[k]}u_\infty)| =  \E_{\infty}(u_\infty)\Big\}\Big)  ,
  \end{split}
\]
which yields
\beq
\label{1.10}
\begin{split}
\ \ \int_{\big\{|\H (\mathrm D^{[k]}u_\infty)| <  \E_{\infty}(u_\infty)\big\}} \big| \H \big(\mathrm D^{[k]}u_\infty \big) \big|^p \, \d \mL^n  \
 =\, \E_{\infty}(u_\infty)^p \mL^n\Big(\Big\{|\H (\mathrm D^{[k]}u_\infty)| < \E_{\infty}(u_\infty)\Big\}\Big).
  \end{split}
\eeq
If we hypothetically had
\[
\mL^n\Big(\Big\{|\H (\mathrm D^{[k]}u_\infty)| < \E_{\infty}(u_\infty)\Big\}\Big) >0,
\]
then identity \eqref{1.10} above results in the contradiction
\[
\begin{split}
\av_{\big\{|\H (\mathrm D^{[k]}u_\infty)| <  \E_{\infty}(u_\infty)\big\}} \big| \H \big(\mathrm D^{[k]}u_\infty \big) \big|^p \, \d \mL^n  \
 =\, \E_{\infty}(u_\infty)^p .
  \end{split}
\]
Therefore, it follows that
\[
\mL^n\Big(\Big\{|\H (\mathrm D^{[k]}u_\infty)| < \E_{\infty}(u_\infty)\Big\}\Big) =0,
\]
which in view of \eqref{1.9} means that $|\H (\mathrm D^{[k]}u_\infty)| = \La$ a.e.\ on $\Om$, with 
\[
 {\La :=  \E_{\infty}(u_\infty) > 0}.
\]
Hence, $u_\infty$ solves \eqref{1.5}, as claimed. This completes the proof.  \qed
\ms

\begin{remark} In the context of Theorem \ref{theorem1} above, it might be reasonable to expect that $\La\geq \La_*$ in part (2), for the same $\La_*$ as in the assumptions. This seems generically speaking to be the case, but it does not appear directly deducible without imposing specific assumptions on $\H$ (instead of relying on the weaker requirement of mere solvability of \eqref{1.5}). Let us demonstrate this in a simple particular case. If for example $|\H|$ is $k$-quasiconvex depending only on the leading term of $k$-th order derivatives, and the boundary condition $\phi$ is a polynomial of degree $k$, then as in Proposition \ref{proposition1} one can define $\La_* $ by \eqref{2.2}. In that case, by H\"older's inequality, the definition of $k$-quasiconvexity, the Gauss-Green theorem and the satisfaction of the Dirichlet boundary condition up to derivatives of order $k-1$, we can indeed estimate
\[
\begin{split}
\La = \underset{\Om}{\ess \, \sup}\,  \big|\H (\mathrm D^{k} u_\infty)\big| &\geq   \av_{\Om} \big|\H (\mathrm D^{k} u_\infty)\big| \, \mathrm d \mL^n
\\
&\geq   \left|\H \! \left(\, \av_{\Om}\mathrm D^{k} u_\infty \, \mathrm d \mL^n\right)\right| 
\\
&=  \left|\H \! \left(\, \av_{\Om}\mathrm D^{k} \phi \, \mathrm d \mL^n\right)\right| 
\\
&= \underset{\Om}{\ess \, \sup}\,  \big|\H (\mathrm D^{k} \phi )\big|
\\
&\geq \La_*.
\end{split}
\]

\end{remark}

\subsection*{Acknowledgement} The author is indebted to Nick Barron (Loyola, Chicago, USA) for posing a few years ago the interesting question discussed in this note, and to Roger Moser (Bath, UK) for scientific discussions relevant to the problem presented herein. He would also like to thank the referee of this paper for the careful reading of the manuscript and their constructive comments, which improved the presentation as well as the content of this work.

\ms



\bibliographystyle{amsplain}

\end{document}